\documentclass[prd, onecolumn,draft]{revtex4}

\usepackage{amsmath}
\usepackage{graphicx}
\usepackage{dcolumn}
\usepackage{bm}
\usepackage{wrapfig}
\usepackage{epsfig}
\usepackage{breqn}
\usepackage{amssymb}
\usepackage{mathtools}
\usepackage{subfigure}
 \usepackage[usenames,dvipsnames]{pstricks}
 \usepackage{epsfig}
 \usepackage{pst-grad} 
 \usepackage{pst-plot} 

\begin{document}


\title{On an Optimal Solution to the Film Scheduling and Showtime Staggering Problem}
\author{Ikjyot Singh Kohli}
	\email{isk@mathstat.yorku.ca}
\affiliation{York University - Department of Mathematics and Statistics}
\author{Katherine Goff Inglis}
\email{kgoff@ryerson.ca}
\affiliation{
Independent Researcher
}

\date{October 22nd, 2019}

\begin{abstract}
The scheduling of films is a major problem for the movie theatre exhibition business. The problem is two-fold: movie exhibitors ideally would like to schedule films to screens in their various locations to maximize attendance and revenue, but would also like to schedule these films such that neighbouring theatre locations play the same films at different times thus giving guests a multitude of showtime options. We refer to this latter problem as the showtime \emph{staggering} problem. We give an exact formulation of this scheduling problem using binary integer linear optimization, and provide a solved example as well. This work further shows that the optimal scheduling of films cannot be done across all theatre locations at once, but rather, must be done for each cluster of neighbouring locations.

\end{abstract}
\maketitle 

\section{Introduction}
In an era of data driven digital transformation, a customer driven business strategy is essential for success.  In the motion picture industry, movie exhibitors must compete to win share of consumers’ entertainment time (and wallet) against digital entertainment alternatives offered by mammoth sized, digital focused, competitors like Netflix, Amazon and Disney \cite{goff2017}.  Customer loyalty, point-of-sale and digital payment platforms produce rich insights that can leveraged to inform business operations and automate the decision-making process, effectively enabling movie exhibitors to compete using analytics and artificial intelligence.  This study presents a new, customer driven, quantitative approach to movie scheduling that can be utilized by movie exhibitors to increase attendance and market share.  

The role of the exhibitor is to show films that are produced by movie studios (see \cite{somlo2005} for more details on the roles of the stakeholders in the movie industry).   Exhibitors do not have decision making authority over the movies that are produced by the studios.  Therefore, one of the most significant business decisions exhibitors make is determining which movies to feature in their multiplexes, also referred to as the movie scheduling problem.  In the literature, the scheduling of movies to screens has been compared to that of scheduling jobs to machines.  However, it is recognized that the movie scheduling problem is non-trivial as it involves several additional constraints not applicable to the standard machine scheduling problem, most notably that the value of the job (movie) is not fixed, it varies based on consumer demand.   As such, sophisticated quantitative methods are required for multiplex movie exhibitors to effectively solve the movie scheduling problem and optimize the scheduling of movies so that demand (and profit) is maximized.  

A number of review papers exist that provide overviews of the various topics and research directions relating to the motion picture industry (e.g., \cite{Moul2005}, \cite{Weinberg2006}, \cite{Eliashbergetal.2008}).  \cite{Goff2} provide a review of scheduling problems and opportunities in motion picture exhibition and use a market driven approach to identify key decision support factors associated with the movie scheduling problem.  This study extends the existing literature by presenting a new quantitative model that addresses the opportunities identified by \cite{Goff2} within the “Network Considerations” and “Show Time Schedule” factors.  Specifically, this paper presents a movie scheduling model that recommends the optimal arrangement of show times for movies and locations, considering the interdependency between theatre locations that are part of the same multiplex chain, located in a common geography.       

The contribution of this work extends existing models in two ways.  Firstly, this is only study that takes an exhibitor’s perspective optimizing a chain of multiplexes (\cite{somlo2011} do this from the distributors perspective).  Secondly, this is the only study that optimizes showtimes across multiple locations in a single geography, providing more options of movie times for customers to choose from in urban markets.  The models presented in the existing literature (e.g. \cite{Eliashberg2001}, \cite{Eliashberg2009a}, \cite{Eliashberg2009b}, \cite{somlo2005}, \cite{somlo2011}) take a micro approach to optimization, that is they optimize profit through scheduling within a single multiplex location.  The model presented in this paper takes a macro approach to optimization: maximizing profit for a market (group of ‘crossover’ locations) rather than individual multiplex locations.  This macro approach is leveraged to generate a variety of showtimes for each movie on each day in a common geography. 


In the next section of the paper, the mathematical model is described. Following that a numerical example is provided.  The paper concludes by identifying next steps for future research to further support movie exhibitors in optimizing their business using data and quantitative methods.

\section{Mathematical Model}
We consider theatre locations that are part of a group of crossover locations, that is, theatre locations that are in neighbouring locations. The problem at hand requires one to schedule films across all screens in each such theatre location such that the scheduling results in a maximization of theatre attendance/revenue. Further, the showtimes are required to be \emph{staggered}. That is, for any given time, no locations can play the same film.

The decision variables that we use in the model are binary decision variables and are denoted: $X_{sf_c} = \{1,0\}$, which denotes that film $f$ is played on screen $s$, with showtime configuration $c$. We further assume that the screens in all neighbouring locations are given a unique identifier. In this paper, we just assume that they are indexed from 1 to the total number of screens in all neighbouring locations which we denote by $S$. The concept of showtime configurations requires some further explanation, which we now give. In our model, based on how movie theatres are known to schedule films to screens, we assume that a given film is assigned to a given screen for a playweek, and throughout a given day, that film plays at several times within the theatre location's operating window. More precisely, let $t_0$ denote the first available time any film can be played in a theatre location, and let $t_f = T$ denote the last possible time a film can be played in the same theatre locations. Then, for any given film, the possible configurations are given by 
\begin{equation}
c = \{t_0, t_0 + \Delta t, t_1 + \Delta t, t_2 + \Delta_t, + \ldots + t_f = T\},
\end{equation}
where $\Delta t$ is a staggering interval, which can be thought of as a time spacing by which we have the flexibility to move showtimes around for a specific film to satisfy the staggering constraint. Clearly, the number of possible configurations is dependent on the film's runtime, which is why we have used the notation $c$ in our decision variables to make this dependency clear. 

The first constraint on the decision variables deals with the scheduling of films to individual screens. For this we require that for $s=1, \ldots S$ screens and $f_c=1, \ldots F$ film configurations.
\begin{eqnarray}
\sum_{f_c=1}^{F} X_{1 f_c} &= 1, \\
\sum_{f_c=1}^{F}  X_{2 f_ c} &= 1, \\
\sum_{f_c=1}^{F}  X_{3 f_c} &= 1, \\
&\vdotswithin{=} \\
\sum_{f_c=1}^{F} X_{S f_c} &=1.
\end{eqnarray}

In addition to this, we have staggering constraints that require a one to one mapping of screens to film configurations. This can be represented by the set of constraints for each film configuration $f_c = 1, \ldots , F$:
\begin{eqnarray}
\sum_{s=1}^{S}  X_{s 1_c} & \leq 1, \\
\sum_{s=1}^{S}  X_{s 2_c} & \leq 1, \\
\sum_{s=1}^{S}  X_{s 3_c} & \leq 1, \\
&\vdotswithin{=} \\
\sum_{s=1}^{S}   X_{s F_ c}  & \leq 1. 
\end{eqnarray}
Note that, in these staggering constraints each different configuration and each different value of $c$ or $f_c$ will introduce a new inequality constraint. The numerical example we employ below makes this clear.

Further, the objective function is given by:
\begin{equation}
Z = \sum_{s=1}^{S} \sum_{f_c=1}^{F} f(s,f_c) X_{sf_c},
\end{equation}
where the coefficients $f(s,f_c)$ are given by a statistical learning model from the data, which represent the predicted attendance for a given film configuration on a given screen in a given theatre location. Building an accurate predictive model, although, not the subject of the current paper, is absolutely key in generating the optimal schedule. In principle, given the amount of ticket sales data that theatres collect, this should not be a difficult thing to do. Further, because our model schedules movies per individual screen, this gives the model much more power in the sense that implicit within the statistical learning function are premium moviegoing experiences such as IMAX. 

\subsection{Numerical Example}
We now present a numerical example to demonstrate the implementation of the optimization problem under consideration. In this example we consider three neighbouring theatre locations, where the first location has screens labeled 1,2,3, the second location has screens labeled 4,5 and the third location has screens labeled 6,7,8,9. Further, we consider the scheduling of 5 films. The films have runtimes of 90 minutes, 85 minutes, 100 minutes, 110 minutes, and 120 minutes respectively. We further assume that the films are being scheduled on a particular day where each theatre location opens at 12:00 PM and has the last possible showtime at 11:00 PM, with a staggering interval of thirty minutes. Therefore, the following showtime configurations are possible:
\begin{enumerate}
	\item Films 1, 2:  (12:30, 2:00, 3:30, 5:00, 6:30, 8:00, 9:30, 11:00), (12:00, 1:30, 3:00, 4:30, 6:00, 7:30, 9:00, 10:30), (1:00, 2:30, 4:00, 5:30, 7:00, 8:30, 10:00)
	\item Films 3, 4, 5: (1:00, 3:00, 5:00, 7:00, 9:00, 11:00), (12:30, 2:30, 4:30, 6:30, 8:30, 10:30), (12:00, 2:00, 4:00, 6:00, 8:00, 10:00), (1:30, 3:30, 7:30, 9:30).
\end{enumerate}
We therefore see that films 1 and 2 have 3 possible time configurations, while films 3, 4, and 5 have 4 possible time configurations.

We also assume that a statistical learning model has been implemented that generates  cumulative attendance predictions for each screen and film configuration combination. There will be 144 such predictions, hence, the objective function will have 144 terms. The objective function would be:
\begin{dmath}
Z = \sum_{s=1}^{S} X_{s  1_1} + X_{s  1_2} + X_{s  2_1} + X_{s  2_2} + X_{s  3_1} + X_{s  3_2} + X_{s  3_3} + X_{s  3_4} 
+ X_{s 4_1} + X_{s  4_2} + X_{s  4_3} + X_{s  4_4} + X_{s  5_1} + X_{s  5_2} + X_{s  5_3} + X_{s  5_4},
\end{dmath}
which we aim to maximize. For the sake of this example, we have provided the full objective function used in this simulation in the appendix.

The constraints are as follows. First, for each individual screen, we have 9 constraint equations:
\begin{dmath}
X_{1  1_1} + X_{1  1_2} + X_{1  2_1} + X_{1  2_2} + X_{1  3_1} + X_{1  3_2} + X_{1  3_3} + X_{1  3_4} 
+ X_{1 4_1} + X_{1  4_2} + X_{1  4_3} + X_{1 4_4} + X_{1  5_1} + X_{1  5_2} + X_{1  5_3} + X_{1  5_4} = 1, 
\end{dmath}
\begin{dmath}
X_{2  1_1} + X_{2  1_2} + X_{2  2_1} + X_{2  2_2} + X_{2  3_1} + X_{2  3_2} + X_{2  3_3} + X_{2  3_4} 
+ X_{2 4_1} + X_{2  4_2} + X_{2  4_3} + X_{2 4_4} + X_{2  5_1} + X_{2  5_2} + X_{2  5_3} + X_{2  5_4} = 1, 
\end{dmath}
\begin{dmath}
X_{3 1_1} + X_{3  1_2} + X_{3  2_1} + X_{3  2_2} + X_{3  3_1} + X_{3  3_2} + X_{3  3_3} + X_{3  3_4} 
+ X_{3 4_1} + X_{3  4_2} + X_{3  4_3} + X_{3 4_4} + X_{3  5_1} + X_{3  5_2} + X_{3  5_3} + X_{3  5_4} = 1, 
\end{dmath}
\begin{dmath}
X_{4 1_1} + X_{4 1_2} + X_{4  2_1} + X_{4  2_2} + X_{4  3_1} + X_{4  3_2} + X_{4  3_3} + X_{4  3_4} 
+ X_{4 4_1} + X_{4 4_2} + X_{4  4_3} + X_{4 4_4} + X_{4  5_1} + X_{4  5_2} + X_{4 5_3} + X_{4  5_4} = 1, 
\end{dmath}
\begin{dmath}
X_{5 1_1} + X_{5 1_2} + X_{5  2_1} + X_{5  2_2} + X_{5  3_1} + X_{5  3_2} + X_{5 3_3} + X_{5  3_4} 
+ X_{5 4_1} + X_{5 4_2} + X_{5  4_3} + X_{5 4_4} + X_{5  5_1} + X_{5  5_2} + X_{5 5_3} + X_{5  5_4} = 1, 
\end{dmath}
\begin{dmath}
X_{6 1_1} + X_{6 1_2} + X_{6  2_1} + X_{6 2_2} + X_{6  3_1} + X_{6  3_2} + X_{6 3_3} + X_{6  3_4} 
+ X_{6 4_1} + X_{6 4_2} + X_{6  4_3} + X_{6 4_4} + X_{6  5_1} + X_{6  5_2} + X_{6 5_3} + X_{6  5_4} = 1, 
\end{dmath}
\begin{dmath}
X_{7 1_1} + X_{7 1_2} + X_{7  2_1} + X_{7 2_2} + X_{7 3_1} + X_{7  3_2} + X_{7 3_3} + X_{7  3_4} 
+ X_{7 4_1} + X_{7 4_2} + X_{7  4_3} + X_{7 4_4} + X_{7  5_1} + X_{7  5_2} + X_{7 5_3} + X_{7 5_4} = 1, 
\end{dmath}
\begin{dmath}
X_{8 1_1} + X_{8 1_2} + X_{8  2_1} + X_{8 2_2} + X_{8 3_1} + X_{8  3_2} + X_{8 3_3} + X_{8  3_4} 
+ X_{8 4_1} + X_{8 4_2} + X_{8 4_3} + X_{8 4_4} + X_{8  5_1} + X_{8  5_2} + X_{8 5_3} + X_{8 5_4} = 1, 
\end{dmath}
\begin{dmath}
X_{9 1_1} + X_{9 1_2} + X_{9  2_1} + X_{9 2_2} + X_{9 3_1} + X_{9  3_2} + X_{9 3_3} + X_{9  3_4} 
+ X_{9 4_1} + X_{9 4_2} + X_{9 4_3} + X_{9 4_4} + X_{9 5_1} + X_{9  5_2} + X_{9 5_3} + X_{9 5_4} = 1.
\end{dmath}

We also implement the staggering constraints by the following set of 16 equations:
where, of course, each $X_{s f_c} \in \{0,1\}$.
\begin{dmath}
X_{1 1_1} + X_{2 1_1} + X_{3 1 _1} + X_{4 1 _1} + X_{5 1 _1}  + X_{6 1 _1}  + X_{7 1 _1}  + X_{8 1 _1}  + X_{9 1 _1}  \leq 1,
\end{dmath}
\begin{dmath}
X_{1 1_2} + X_{2 1_2} + X_{3 1_2} + X_{4 1_2} + X_{5 1_2} + X_{6 1_2} + X_{7 1_2}  + X_{8 1_2}  + X_{9 1_2} \leq 1, 
\end{dmath}
\begin{dmath}
X_{1 2_1} + X_{2 2_1} + X_{3 2_1} + X_{4 2_1} + X_{5 2_1} + X_{6 2_1} + X_{7 2_1} + X_{8 2_1} + X_{9 2_1} \leq 1,
\end{dmath}
\begin{dmath}
X_{1 2_2} + X_{2 2_2} + X_{3 2_2} + X_{4 2_2} + X_{5 2_2} + X_{6 2_2} + X_{7 2_2}  + X_{8 2_2}  + X_{9 2_2}  \leq 1,
\end{dmath}
\begin{dmath}
X_{1 3_1} + X_{2 3_1} + X_{3 3_1} + X_{4 3_1} + X_{5 3_1} + X_{6 3_1} + X_{7 3_1} + X_{8 3_1} + X_{9 3_1} \leq 1,
\end{dmath}
\begin{dmath}
X_{1 3_2} + X_{2 3_2} + X_{3 3_2} + X_{4 3_2} + X_{5 3_2} + X_{6 3_2} + X_{7 3_2} + X_{8 3_2} + X_{9 3_2} \leq 1,
\end{dmath}
\begin{dmath}
X_{1 3_3} + X_{2 3_3}  + X_{3 3_3}  + X_{4 3_3}  + X_{5 3_3}  + X_{6 3_3}  + X_{7 3_3} + X_{8 3_3}  + X_{9 3_3}  \leq 1,
\end{dmath}
\begin{dmath}
X_{1 3_4} + X_{2 3_4} + X_{3 3_4} + X_{4 3_4} + X_{5 3_4} + X_{6 3_4} + X_{7 3_4} + X_{8 3_4} + X_{9 3_4} \leq 1,
\end{dmath}
\begin{dmath}
X_{1 4_1} + X_{2 4_1}  + X_{3 4_1} + X_{4 4_1} + X_{5 4_1} + X_{6 4_1}  + X_{7 4_1}  + X_{8 4_1}  + X_{9 4_1}  \leq 1,
\end{dmath}
\begin{dmath}
X_{1 4_2} + X_{2 4_2} + X_{3 4_2}  + X_{4 4_2}  + X_{5 4_2}  + X_{6 4_2}  + X_{7 4_2}  + X_{8 4_2}  + X_{9 4_2}  \leq 1,
\end{dmath}
\begin{dmath}
X_{1 4_3} + X_{2 4_3} + X_{3 4_3} + X_{4 4_3} + X_{5 4_3} + X_{6 4_3} + X_{7 4_3} + X_{8 4_3} + X_{9 4_3} \leq 1,
\end{dmath}
\begin{dmath}
X_{1 4_4} + X_{2 4_4}  + X_{3 4_4}  + X_{4 4_4}  + X_{5 4_4}  + X_{6 4_4} + X_{7 4_4}  + X_{8 4_4}  + X_{9 4_4}  \leq 1,
\end{dmath}
\begin{dmath}
X_{1 5_1} + X_{2 5_1} + X_{3 5_1}  + X_{4 5_1}  + X_{5 5_1} +  X_{6 5_1}   + X_{7 5_1}  + X_{8 5_1}  + X_{9 5_1}  \leq 1,
\end{dmath}
\begin{dmath}
X_{1 5_2} + X_{2 5_2} + X_{3 5_2}  + X_{4 5_2}  + X_{5 5_2} +  X_{6 5_2}   + X_{7 5_2}  + X_{8 5_2}  + X_{9 5_2}  \leq 1,
\end{dmath}
\begin{dmath}
X_{1 5_3} + X_{2 5_3} + X_{3 5_3}  + X_{4 5_3}  + X_{5 5_3} +  X_{6 5_3}   + X_{7 5_3}  + X_{8 5_3}  + X_{9 5_3}  \leq 1,
\end{dmath}
\begin{dmath}
X_{1 5_4} + X_{2 5_4} + X_{3 5_4}  + X_{4 5_4}  + X_{5 5_4} +  X_{6 5_4}   + X_{7 5_4}  + X_{8 5_4}  + X_{9 5_4}  \leq 1.
\end{dmath}

Numerically solving the maximization problem (13)-(38), we obtain the solution:
\begin{equation}
X_{1 5_4} = 1, \quad X_{2 5_1} = 1, \quad X_{3 3_2} = 1, \quad X_{4 3_4} = 1, \quad X_{5 2_1} = 1, \quad X_{6 1_2} = 1, \quad X_{7 3_1} = 1, \quad X_{8 5_2} = 1, \quad X_{9 4_4} = 1.
\end{equation}
For clarity, we present the result in the table below:
\begin{table}[h]
\begin{center}\begin{tabular}{|c|c|c|}\hline Screen & Film & Configuration \\\hline 1 & 5 & 4 \\\hline 2 & 5 & 1 \\\hline 3 & 3 & 2 \\\hline 4 & 3 & 4 \\\hline 5 & 2 & 1 \\\hline 6 & 1 & 2 \\\hline 7 & 3 & 1 \\\hline 8 & 5 & 2 \\\hline 9 & 4 & 4 \\\hline 
\end{tabular} 
\caption{Solution of Simulated Optimization Staggering Problem}
\end{center}
\end{table}
One can see from Table I, that the five films have been allocated to the 9 different screens in the various theatre locations in a way that maximizes attendance, but also in a way that no films have the same configuration schedule. For example, screens 1, 2, and 8 are all scheduled to play film 5, but with time configurations 4,1, and 2. Similarly, screens 3,4, and 7 are all scheduled to play film 3, but with time configurations 2,4, and 1. Hence, the film showtimes are staggered, which was one of the goals of this analysis.

\section{Conclusions}
The mathematical model presented in this paper provides a novel approach to scheduling films and showtimes in multiplexes that are part of a cinema chain.  The implementation of this approach in practice would benefit both movie exhibitors and consumers.  For movie exhibitors, this algorithm enables managers to monetize data by applying quantitative methods to generate movie schedules that yield the highest attendance numbers.  Implementation of this approach could impact exhibitor profit in two ways.  Firstly, this approach would lead to increased revenues as a result of incremental ticket sales.  Secondly, this method could result in decreased labour costs by utilizing a more efficient process of schedule generation.  Consumers would benefit from having this algorithm implemented in practice, they would have more choices of showtimes for the movies they want to see and more relevant movies scheduled at the theatres they frequent.  To quantify the potential impact of this work it is necessary to pilot the algorithm in a real-world scenario with a movie exhibitor, this is the most important next step we have identified for this research.  In addition to this, there are several other research opportunities to further extend our model that are discussed in following paragraph.

Movie exhibitors generate revenues through three main avenues: ticket sales, concession food sales and advertising sales.  People go to movie theatres to see particular movies, therefore concession food and advertising sales are ancillary and dependent on the movies that are playing.  There are opportunities to extend our model and develop algorithms to deliver the best outcome for food and advertising sales.  For food sales, the desired outcome is similar to ticket sales: to yield the maximum amount of food profit for a neighbourhood. 

Our model can be extended to account for this by adjusting the input box office demand forecast to a combined profit forecast including ticket and food sales.  This model extension would account for the difference in food offerings between neighbouring cinema locations and within cinema auditoriums (e.g. Alcohol being served in VIP).  It was not within the scope of this paper to produce a box office demand forecast but is relevant to note that this is another fruitful area for future research.  With respect to advertising sales there is a unique opportunity to extend our model to generate (and optimize) the scheduling of advertisements and movie trailers that are featured prior to the movie.  The existing literature on this topic is limited.  Some of the existing models incorporate a constant time (e.g. 15 minutes before showtime) to be allocated to trailers and advertisements.  We did not find any papers that generate pre-show schedules that consider the variability in volume and length of movie trailers and advertisements.  An algorithm such as this could facilitate the scheduling and delivery of advertisements based on industry standard CPM or impression-based sales models.  The development and implementation of a scheduling application that considers all profit avenues (box office, concession and advertising) and recommends the optimal outcome for a multiplex chain is an example of how movie exhibitors can compete in a digital era using artificial intelligence and data driven quantitative methods.

\section{Appendix}
The objective function we maximized in the numerical example is given by:
\begin{dmath}
Z= 226 \text{X111}+245 \text{X112}+232 \text{X121}+256 \text{X122}+202 \text{X131}+208 \text{X132}+272 \text{X133}+244 \text{X134}+238 \text{X141}+222 \text{X142}+212 \text{X143}+222 \text{X144}+209 \text{X151}+247 \text{X152}+233 \text{X153}+286 \text{X154}+213 \text{X211}+276 \text{X212}+211 \text{X221}+229 \text{X222}+206 \text{X231}+251 \text{X232}+260 \text{X233}+205 \text{X234}+217 \text{X241}+297 \text{X242}+202 \text{X243}+266 \text{X244}+298 \text{X251}+259 \text{X252}+244 \text{X253}+221 \text{X254}+205 \text{X311}+288 \text{X312}+252 \text{X321}+244 \text{X322}+236 \text{X331}+293 \text{X332}+266 \text{X333}+218 \text{X334}+288 \text{X341}+261 \text{X342}+245 \text{X343}+247 \text{X344}+247 \text{X351}+268 \text{X352}+219 \text{X353}+260 \text{X354}+253 \text{X411}+247 \text{X412}+266 \text{X421}+268 \text{X422}+234 \text{X431}+249 \text{X432}+211 \text{X433}+282 \text{X434}+263 \text{X441}+266 \text{X442}+268 \text{X443}+201 \text{X444}+214 \text{X451}+261 \text{X452}+220 \text{X453}+216 \text{X454}+234 \text{X511}+256 \text{X512}+291 \text{X521}+218 \text{X522}+229 \text{X531}+248 \text{X532}+234 \text{X533}+251 \text{X534}+274 \text{X541}+242 \text{X542}+207 \text{X543}+209 \text{X544}+228 \text{X551}+246 \text{X552}+289 \text{X553}+211 \text{X554}+274 \text{X611}+293 \text{X612}+223 \text{X621}+263 \text{X622}+284 \text{X631}+211 \text{X632}+284 \text{X633}+228 \text{X634}+278 \text{X641}+289 \text{X642}+255 \text{X643}+273 \text{X644}+237 \text{X651}+252 \text{X652}+218 \text{X653}+200 \text{X654}+273 \text{X711}+213 \text{X712}+272 \text{X721}+203 \text{X722}+292 \text{X731}+242 \text{X732}+260 \text{X733}+252 \text{X734}+283 \text{X741}+263 \text{X742}+207 \text{X743}+252 \text{X744}+208 \text{X751}+259 \text{X752}+219 \text{X753}+280 \text{X754}+254 \text{X811}+271 \text{X812}+265 \text{X821}+214 \text{X822}+200 \text{X831}+284 \text{X832}+251 \text{X833}+274 \text{X834}+269 \text{X841}+215 \text{X842}+229 \text{X843}+218 \text{X844}+224 \text{X851}+285 \text{X852}+248 \text{X853}+217 \text{X854}+293 \text{X911}+273 \text{X912}+258 \text{X921}+205 \text{X922}+245 \text{X931}+264 \text{X932}+261 \text{X933}+219 \text{X934}+270 \text{X941}+286 \text{X942}+228 \text{X943}+295 \text{X944}+256 \text{X951}+257 \text{X952}+200 \text{X953}+267 \text{X954}.
\end{dmath}

\newpage 
\bibliography{sources}

\end{document}